\documentclass[%
 reprint,
showpacs, showkeys,
 amsmath,amssymb,
 aps,
]{revtex4-1}

\usepackage{graphicx}
\usepackage{dcolumn}
\usepackage{bm}

\newtheorem {thm}{Theorem}[section]
\newtheorem {lemma}[thm]{Lemma}

\newcommand{\bbmatrix}[1]{\left[ \begin{array}{cccccccccccccccccc} #1 \end{array} \right]}
\newcommand{\SPEC}{{\operatorname{spec}}}

\def\squarebox#1{\hbox to #1{\hfill\vbox to #1{\vfill}}}
\newcommand{\qed}{\hspace*{\fill}
\vbox{\hrule\hbox{\vrule\squarebox{.667em}\vrule}\hrule}\smallskip}

\begin{document}

\preprint{APS/123-QED}

\title{A new method for solving completely integrable PDEs}

\author{Andrey Melnikov}
 \email{andreymath@gmail.com}
\affiliation{
	Department of Mathematics, \\
	33rd and Market Streets \\
	Philadelphia PA, 19104 USA
}

\begin{abstract}
The inverse scattering theory is a basic tool to solve linear differential equations and some Partial Differential Equations (PDEs).
Using this theory the Korteweg-de Vries (KdV), the family of evolutionary Non Linear Schr\" odinger (NLS) equations, 
Kadomtzev-Petviashvili and many more completely integrable PDEs of mathematical physics are solved,
using Zacharv-Shabath scheme.  This last approach includes the use of a Lax pair, and has an advantage to be applied to wider class of equations, like difference equations, but
has a disadvantage to be used only for ``rapidly decreasing solutions''. This technique is also intimately related to completely integrable systems. 
The identifying process of a Lax pair, a system and finally the ``scattering data'' is usually a difficult
process, simplified in many cases by physicals models providing clues of how the scattering data should be chosen.
In this work we show that the scattering data can be encoded into singularities of a very special mathematical object: $J$-unitary, identity at infinity matrix-valued function, 
which, if evolved once, solves the inverse scattering theory, 
if evolved twice solves evolutionary PDEs. The provided scheme seems to be universal in the sense that many (if not all) completely
integrable PDEs arise (or should arise) in this manner by changing the so called ``vessel parameters'' (for examples, solutions of KdV and evolutionary NLS equations are presented). 
The results presented here allow to study different flows (commuting and non-commuting) in a unified approach and provide a rich mathematical arsenal to study
these equations. The results are easily generalized to completely integrable PDEs of $n$ ($\leq\infty$) variables.
\end{abstract}


\maketitle


\section{\label{sec:Intro}Introduction.}
A standard technique to solve some non linear PDEs, initiated in \cite{bib:GGKM} and generalized to a Zacharov-Shabbath scheme in \cite{bib:ZakSha74short} includes
for example the Korteweg-de Vries equation
\begin{equation} \label{eq:KdV}
q_t = - \dfrac{3}{2} q q_x + \dfrac{1}{4} q_{xxx},
\end{equation}
where $q_t, q_x$ denote the partial derivatives. This equation arises in the study of waves on shallow water surfaces.
A second important example is the Non Linear Schr\" odinger (NLS) evolutionary equation:
\begin{equation} \label{eq:ENLS}
i \beta_t + \beta_{xx} + 2 |\beta|^2 \beta = 0, \quad \beta(x,0) = \beta(x),
\end{equation}
where $\beta=\beta(x,t)$ is a complex valued function of two real variables $x,t$ and $\beta(x,0)=\beta(x)$ is the initial condition, defined usually on $\mathbb R$.
This equation plays a special role in optics and water waves. we show how these two important equations of mathematical physics arise as special cases of a theory of vessels,
recently developed by the author, initiated by works of M. Liv\v sic \cite{bib:Vortices}. As we would like to focus on the technique and not on mathematical aspects
of the theory, we will not actually define the vessel itself, but will use all its necessary ingredients to construct solutions of many completely integrable PDEs.

\section{\label{sec:describ}Description of the method}
\subsection{Scattering theory}
Suppose that $p\in\mathbb N$ and let
\[ \sigma_1=\sigma_1^*, \quad \sigma_2=\sigma_2^*, \quad \gamma=-\gamma^*,
\]
be three $p\times p$ matrices such that $\sigma_1$ is invertible. We call such a triple as \textbf{vessel parameters}. 
For example, some triples for the case $p=2$ are given in Table \ref{tab:3BasicEx}. 
\begin{table}[b]
\caption{\label{tab:3BasicEx}Parameters corresponding to Sturm-Liouville (SL), Non Linear Schr\" odinger (NLS)
Equation and Canonical systems. $I$ - the identity matrix, $0$ - the zero matrix.}
\begin{tabular}{l|c|c|c|}
    p=2 & $\sigma_1$& $\sigma_2$ & $\gamma$ \\
\colrule
SL & $\bbmatrix{0&1\\1&0}$& $\bbmatrix{1&0\\0&0}$& $\bbmatrix{0&0\\0&i}$ \\
NLS & $I$ & $\bbmatrix{\frac{1}{2}&0\\0&-\frac{1}{2}}$ & $0$ \\
Can. Sys. & $\bbmatrix{0&i\\-i&0}$ & $I$ & $0$ \\
\end{tabular}
\end{table}
By a \textit{scattering data} we mean a $p\times p$ function $S(\lambda)$, called
\textbf{scattering data matrix} possessing a realization
\begin{multline} \label{eq:S0realized}
S(\lambda) = I - B^*_0 \mathbb X_0^{-1}(\lambda I - A)^{-1}B_0\sigma_1, \\
A \mathbb X_0 + \mathbb X_0 A^* + B_0\sigma_1 B_0 = 0, \quad \mathbb X_0^* = \mathbb X_0
\end{multline}
where for an auxiliary Hilbert space $\mathcal H$ the operators act as follows: $\mathbb X_0, A: \mathcal H\rightarrow\mathcal H$, $B_0:\mathbb C^p\rightarrow\mathcal H$.
Such realizations where studied in \cite{bib:bgr} and we present the \textit{regular case, where all the operators are bounded} to clarify explanations and to simplify formulas,
but the general case with an unbounded $A$
is presented later in Section \ref{sec:NonReg}. Letting the operators evolve with $x$ as follows:
\begin{eqnarray}
\label{eq:DB} \frac{d}{dx} B(x)\sigma_1 = - A B(x)\sigma_2-B(x)\gamma, \quad B(0)=B_0, \\
\label{eq:DX} \frac{d}{dx} \mathbb X(x) = B(x) \sigma_2 B^*(x), \quad \mathbb X(0)=\mathbb X_0
\end{eqnarray}
we obtain that for the self-adjoint $\mathbb X(x)$ the Lyapunov equation
\begin{equation} \label{eq:Lyapunov}
A \mathbb X(x) + \mathbb X(x) A^* + B(x)\sigma_1 B(x) = 0
\end{equation}
holds (see \cite{bib:SchurIEOT}) and the function
\begin{equation} \label{eq:Srealized}
S(\lambda,x) = I - B^*(x) \mathbb X^{-1}(x)(\lambda I - A)^{-1}B(x)\sigma_1
\end{equation}
serves as a B\" acklund transformation:
\begin{thm}[\cite{bib:SLVessels, bib:SchurIEOT}]\label{thm:Backlund}
Fix $\lambda\not\in\SPEC(A)$ and let $u(\lambda,x)$ be a solution of the input Linear Differential Equation (LDE) 
\begin{equation} \label{eq:InCC}
		-\sigma_1\dfrac{\partial}{\partial x}u(\lambda,x) + (\sigma_2 \lambda + \gamma)u(\lambda,x) = 0.
\end{equation}
Then the function $y(\lambda,x)=S(\lambda,x)u(\lambda,x)$ is differentiable and satisfies the output LDE
\begin{equation} \label{eq:OutCC}
		-\sigma_1\dfrac{\partial}{\partial x} y(\lambda,x) + (\sigma_2 \lambda + \gamma_*(x))y(\lambda,x) = 0,
\end{equation}
where using $H_0(x) = B^*(x)\mathbb X^{-1}(x) B(x)$
\begin{equation} \label{eq:Linkage}
\gamma_*(x) = \gamma + \sigma_2 H_0(x) \sigma_1 - \sigma_1 H_0(x) \sigma_2.
\end{equation}
\end{thm}
In order to explain this Theorem, the Table \ref{tab:3BasicEx} is useful. Taking the vessel parameters
$\sigma_1, \sigma_2, \gamma$ from a row in this Table and using an arbitrary function $S(\lambda)$ as in \eqref{eq:S0realized},
we will obtain that $S(\lambda,x)$, defined in \eqref{eq:Srealized} maps solutions of a simple LDE, with constant coefficients,
to a more complicated one, defined by $\gamma_*(x)$. This is an analogue of the scattering operator. Moreover, in SL case
it is also possible to characterize precisely the equation \eqref{eq:OutCC}. It turns out \cite{bib:SLVessels} that 
for $y(\lambda,x)=\bbmatrix{y_1(\lambda,x)\\y_2(\lambda,x)}$ the second entry
$y_2(\lambda,x)$ is uniquely derived from $y_1(\lambda,x)$, and the first one satisfies the Sturm-Liouville LDE
\[ -\frac{\partial^2}{\partial x^2} y_1 + q(x) y_1(x) = -i \lambda y_1,
\]
where $q(x) = -2\frac{d^2}{dx^2} \ln (\det \mathbb X_0^{-1}\mathbb X(x))$ (see \cite{bib:SLVessels} for details). Worth noticing
that the famous \textbf{tau}-function is defined in general on the basis of this example as follows:
\begin{equation} \label{eq:tau} \tau(x) = \det(\mathbb X_0^{-1}\mathbb X(x))
\end{equation}

\subsection{\label{sec:SLKdV}Evolutionary equations of SL-KdV type}
An interesting phenomenon occurs when we let the operators evolve with respect to a variable $t$. We choose to represnt
the following two evolutions:
\begin{eqnarray} 
\label{eq:DBt}
 \dfrac{\partial}{\partial t} B(x,t) = i A \dfrac{\partial}{\partial x} B(x,t), \quad B(x,0)=B(x), \\
\label{eq:DXt} \dfrac{\partial}{\partial t} \mathbb X = 
 i A B \sigma_2 B^* - i B\sigma_2 B^* A^* +
 i B\gamma B^*, \\
 \nonumber \mathbb X(x,0)=\mathbb X(x),
\end{eqnarray}
which in the case of SL parameters create solutions of the KdV equation \eqref{eq:KdV}. 
That's why we also call these evolutionary equations as \textit{of SL-KdV type}.
In the case when all the operators are bounded it is a matter of simple calculations to check
the following
\begin{thm}[\cite{bib:KdVVessels}] For the SL vessel parameters, let $B(x,t), \mathbb X(x,t)$ be two operators satisfying
\eqref{eq:DB}, \eqref{eq:DX}, \eqref{eq:DBt}, \eqref{eq:DXt}. Then the function
$q(x,t) = -2 \frac{d^2}{dx^2} \ln (\det \mathbb X_0^{-1} \mathbb X(x))$ satisfies \eqref{eq:KdV}.
\end{thm}
\begin{thm}[\cite{bib:ENLS}] For the NLS vessel parameters, 
let $B(x,t), \mathbb X(x,t)$ be two operators satisfying
\eqref{eq:DB}, \eqref{eq:DX}, \eqref{eq:DBt}, \eqref{eq:DXt}. Then the function
$\beta(x,t) = \bbmatrix{0&1} \gamma_*(x,t)\bbmatrix{1\\0}$ satisfies \eqref{eq:ENLS} for
$\gamma_*(x,t)$, defined by \eqref{eq:Linkage}.
\end{thm}
If the reader is interested in what happens for the canonical systems parameters, 
the answer is as follows and its proof is omitted. In order to understand how it is proved, see Section 
\ref{sec:evolution}.
\begin{thm} For the canonical systems vessel parameters, 
let $B(x,t), \mathbb X(x,t)$ be two operators satisfying
\eqref{eq:DB}, \eqref{eq:DX}, \eqref{eq:DBt}, \eqref{eq:DXt}. Then 
$\gamma_*(x,t)=\bbmatrix{-2 i \beta(x,t)&i h(x,t)\\ih(x,t)&2i\beta(x,t)}$ with $h(x,t)=\sqrt{\dfrac{1}{(x+K)^2} - 4\beta^2}$ ($K\in\mathbb R$)
and the real-valued function $\beta(x,t)$ satisfies:
\begin{equation} \label{eq:CanSysE}
\dfrac{\partial}{\partial t}\big( \sqrt{\dfrac{1}{(x+K)^2} - 4\beta^2}\big) = -2 \dfrac{\beta}{(x+K)^2} + \beta_{xx},
\end{equation}
\end{thm}
Finally, the class of potentials, corresponding to the regular vessels is as follows:
\begin{thm}\label{thm:AnalPot}
Suppose that $B(x), \mathbb X(x)$ are defined by \eqref{eq:DB}, \eqref{eq:DX} using a bounded operator $A$, then
$\gamma_*(x)$ is infinitely differentiable at each point, where $\tau(x)\neq 0$.
\end{thm}
\textbf{Proof:} Notice that using \eqref{eq:DB}, \eqref{eq:DX} the expression \eqref{eq:Linkage} is infinitely times
differentiable.
\qed

\subsection{Commuting flows. Hierarchies}
The formulas \eqref{eq:DBt}, \eqref{eq:DXt} have a special property. It is almost readable that the mixed second derivatives,
applied to $B(x,t)$ are equal: $\frac{\partial^2}{\partial x\partial y}B(x,t) = \frac{\partial^2}{\partial y\partial x}B(x,t)$. But the same
is true for $\mathbb X(x,t)$, which is shown in \cite{bib:KdVHierarchy}:
\begin{lemma} Suppose that $B(x,t),\mathbb X(x,t)$ satisfy \eqref{eq:DB}, \eqref{eq:DX}, \eqref{eq:DBt}, \eqref{eq:DXt}. Then
\[ \frac{\partial^2}{\partial x\partial y} \mathbb X(x,t) = \frac{\partial^2}{\partial y\partial x}\mathbb X(x,t).
\]
\end{lemma}
Similar formulas may be used to derive evolutionary PDEs, constructed from other commuting flows. Let us demand
\begin{equation} \label{eq:nDBt}
 \dfrac{\partial}{\partial t}(B\sigma_1) = A \sum\limits_{i=0}^n A^i B m_i , \quad \bar m_n = (-1)^n m_n 
\end{equation}
and
\begin{equation} \label{eq:nDXt}
\mathbb X'_t = - \sum\limits_{i=0}^n \mathbb Y_i.
\end{equation}
where
\begin{equation}
\label{eq:defY} \mathbb Y_n = \sum\limits_{i=0}^n (-1)^i A^{n-i} B m_n B^* (A^*)^i,
\end{equation}
It was obtained in \cite[Lemma 14]{bib:KdVHierarchy} that in a basic case the corresponding flows commute:
\begin{thm} Suppose that $B(x,t)$ satisfies \eqref{eq:DB} and
\begin{equation} \label{eq:BtKdVn}
 B'_t = (\sqrt{-1}A)^n B'_x.
\end{equation}
Suppose also that $\mathbb X(x,t)$ satisfies \eqref{eq:DX} and \eqref{eq:nDXt}, then
\[
 \frac{\partial^2}{\partial x\partial y}B = \frac{\partial^2}{\partial y\partial x}, \quad
\frac{\partial^2}{\partial x\partial y} \mathbb X = \frac{\partial^2}{\partial y\partial x}\mathbb X.
\]
\end{thm}
The same result conjecturally holds for \eqref{eq:nDBt}, \eqref{eq:nDXt}, when the polynomial with matrix coefficients
$p_0(\lambda) = \lambda\sigma_2+\gamma$ commutes with $\sum\limits_{i=0}^n \lambda^i m_i$.
In a particular case of commuting flows it was also discovered \cite[Theorem 19]{bib:KdVHierarchy} the corresponding family of PDEs,
which also constitute a part of the KdV hierarchy:
\begin{thm} Suppose that $B(x,t), \mathbb X(x,t)$ satisfy \eqref{eq:DB}, \eqref{eq:BtKdVn}, \eqref{eq:DX} and \eqref{eq:nDXt}.
Assume that $b_0 = -\dfrac{1}{4} \beta_{xxx} + \dfrac{3}{2} (\beta_x)^2$ and define recursively
a differential polynomial in $\beta(x)$ as follows ($n=0,1,2\ldots$)
\begin{equation} \label{eq:dbn+1}
4 (b_{n+1})_x = -i (b_n)_{xxx} + 4i (\beta_x b_n)_x.
\end{equation}
Then the 1,2 entry of $\gamma_*(x,t)$, defined by \eqref{eq:Linkage} satisfies
\[ \beta_t = (b_n)_x,
\]
coinciding for $n=0$ with the (integration with respect to $x$ of) the KdV equation \eqref{eq:KdV}.
\end{thm}
\subsection{Completely integrable systems corresponding to SL-KdV type}
Let $(\mathtt{u}(x,t)/\mathtt{x}(x,t)/\mathtt{y}(x,t))$ 
be the (input/ state/ output) triple of the following system of equations
\begin{equation} \label{CISystem}
\left\{ \begin{array}{lll}
\frac{\partial}{\partial t} {\mathtt x}(t,x) = A {\mathtt x}(t,x) + B(x,t) \sigma_1 {\mathtt u}(t,x), \\
\frac{\partial}{\partial x} {\mathtt x}(t,x) = B(x,t) \sigma_2 {\mathtt u}(t,x), \\
{\mathtt y}(t,x) = {\mathtt u}(t,x) - B^*(x,t) \mathbb X^{-1}(x,t) {\mathtt x}(t,x).
\end{array} \right. \end{equation}
which is \textbf{overdetermined}, namely only $\mathtt{u}(t,x)$ satisfying
\[ [- \sigma_1 \frac{\partial}{\partial x} + \sigma_2 \frac{\partial}{\partial t} + \gamma]\mathtt{u}(x,t)  = 0 \]
is allowed. Then the output satisfies
\[
[- \sigma_1 \frac{\partial}{\partial x} + \sigma_2 \frac{\partial}{\partial t} + \gamma_*(x,t)] \mathtt{y} (x,t)  = 0, \]
which can be checked using some algebraic manipulations analogously to the proof of Theorem \ref{thm:Backlund} (B\" acklund transformation).
This system is integrable (for legal inputs!), because the mixed second derivatives are equal:
\[ \frac{\partial^2}{\partial x\partial t} x(t,x) = \frac{\partial^2}{\partial t\partial x} x(t,x), \]
which can be checked by the use of system equations \eqref{CISystem}, the equations \eqref{eq:DB}, \eqref{eq:DX} and \eqref{eq:Lyapunov}
(see \cite{bib:LKMV, bib:Vortices} for details) .
Performing here separation of variables $\mathtt u(t,x) = e^{t \lambda} u(\lambda, x)$,
$\mathtt x(t,x) = e^{t \lambda} x(\lambda, x)$,
$\mathtt y(t,x)= e^{t \lambda} y(\lambda, x)$
we will arrive to the differential equations \eqref{eq:InCC}, \eqref{eq:OutCC} for $u(\lambda, x)$, $y(\lambda, x)$ respectively.
So, vessels indeed come from completely integrable systems by separation of variables.

\subsection{\label{sec:evolution}Evolution of $S(\lambda,x)$ for SL-KdV type}
Notice that the B\" acklund transformation Theorem \ref{thm:Backlund} is equivalent to the fact that
$S$ satisfies
\begin{equation} \label{eq:DS}
\frac{\partial}{\partial x} S = \sigma_1^{-1} (\sigma_2 \lambda + \gamma_*) S - S \sigma_1^{-1} (\sigma_2 \lambda + \gamma)
\end{equation}
Moreover, for the equations of SL-KdV type it is possible to show \cite[(33)]{bib:KdVHierarchy} that
\begin{equation} \label{eq:DSt}
\dfrac{\partial}{\partial t} S = i\lambda \dfrac{\partial}{\partial x} S +
i \dfrac{\partial}{\partial x} [H_0] \sigma_1 S,
\end{equation}
where $H_0$ is the zero moment, while a general moment is defined as $H_n=B^*\mathbb X^{-1}A^n B$. Finally calculating
\[ \dfrac{\partial}{\partial t} \gamma* = \dfrac{\partial}{\partial t} \sigma_2 B^*\mathbb X^{-1}B\sigma_1 -  \dfrac{\partial}{\partial t} \sigma_1 B^*\mathbb X^{-1}B\sigma_2
\]
using \eqref{eq:DBt}, \eqref{eq:DXt} and then rearranging the terms, we will obtain that
\begin{equation} \label{eq:Dgamma*tKdV}
(\gamma_*)'_t = - i \gamma_* (H_0)'_x\sigma_1 + i \sigma_1 (H_0)''_{xx} \sigma_1 +i \sigma_1 (H_0)'_x \gamma_*
\end{equation}
It is a matter of a little effort (after studying the moments) to show that \eqref{eq:Dgamma*tKdV} is equivalent to \eqref{eq:KdV} in the SL case,
to \eqref{eq:ENLS} in for the NLS vessel parameters, and to \eqref{eq:CanSysE} for the canonical system parameters.
\subsection{\label{sec:nDsystems} Higher dimensional completely integrable PDEs}
It is a very interesting question which PDEs are solved, if we choose to use the family of NLS equations,
defined by A.P. {F}ordy, P.P. {K}ulish \cite{bib:NLS}. These also constitute examples of applicability of the suggested
scheme to solve the scattering theory for different values of $p$.

Another important example is the Kadomtzev-Petviashvili equation \cite{bib:KadPetv}. It is a generalization of the KdV equation \eqref{eq:KdV}.
An example of construction of solutions for this equation is as follows. Define for $B=B(x,t,y)$
\[ B_t = B_{xxx} + B_x, \quad \mathbb X_t = \text{\eqref{eq:DXt}} + \text{\eqref{eq:DX}}
\]
and $B_y = B_x, \mathbb X_y = \mathbb X_x$. Then it follows using ideas in \cite{bib:KdVVessels} that
\[ (\gamma_*)_{yy} = \dfrac{\partial}{\partial x} \big[\gamma_t + i \gamma_* (H_0)'_x\sigma_1 - i \sigma_1 (H_0)''_{xx} \sigma_1 - i \sigma_1 (H_0)'_x \gamma_*\big]
\]
which turns to be equal to $(\gamma_*)_{xx}$, resulting in KP equation, when all the terms are transfered to the left hand side.
\section{\label{sec:NonReg}The non regular case}
The same formulas can be used when we use an unbounded operator, giving us a richer family of ``generalized potentials''
$\gamma_*$. In this case the operator $A:D(A)\rightarrow\mathcal H$ is assumed to have a dense domain $D(A)\subseteq\mathcal H$ and constitute a generator of a $C_0$
semi-group. Let us denote its resolvent by $R(\lambda) = (\lambda I - A)^{-1}$, then for $B(x)$ we demand
\begin{equation} \label{eq:ResDB}
0  =  \frac{\partial}{\partial x} [R(\lambda) B(x)] \sigma_1 + A R(\lambda) B(x) \sigma_2  + R(\lambda) B(x).
\end{equation}
Notice that we obtain \textbf{regularity assumptions}, necessary for the equation \eqref{eq:ResDB} to be well defined:
\[ R (\lambda) B(x) \sigma_2 \in D (A) , \quad R (\lambda) B(x) \gamma \in \mathcal H,
\]
holding for all $\lambda\not\in\SPEC(A)$.
The equation \eqref{eq:DX} remains the same and the Lyapunov equation \eqref{eq:Lyapunov} becomes
\begin{multline} \label{eq:ResLyapunov}
\mathbb X(x) R^*(-\bar\lambda) + R(\lambda) \mathbb X(x) - \\
 R(\lambda) B(x) \sigma_1 B^*(x) R^*(-\bar\lambda)  = 0.
\end{multline}
Notice that \eqref{eq:ResDB}, \eqref{eq:ResLyapunov} coincide with \eqref{eq:DB}, \eqref{eq:DX} after cancellations when $A$ is bounded.
In this case, it is possible to show that $S(\lambda,x)$ defined in \eqref{eq:Srealized} will map solutions of \eqref{eq:InCC}
to \eqref{eq:OutCC}, defined by $\gamma_*$, which usually fails to be differentiable, generalizing Theorem \ref{thm:AnalPot}.

\section{Standard models}
\paragraph{Solitons}
The solitons are obtained if we consider so called finite dimensional realizations ($\dim\mathcal H<\infty$). In these cases
the scattering data matrix $S(\lambda)$ is rational. For the Sturm-Liouville case see \cite{bib:SLVessels}.
\paragraph{Spectrum on a curve} If we fix a symmetric with respect to the imaginary axis Jordan curve $\Gamma$
and consider
\[ \mathcal H=L^2(\Gamma) = \{ f(\mu) \mid \int_\Gamma |f(\mu)|^2 < \infty\}.
\]
and $A = i\mu$ - multiplication by the variable of $\Gamma$.
Then there are explicit formulas for the operators $B(x) = B(\mu,x)$, $\mathbb X(x)$ in \cite{bib:GenVessel}, showing
that $S(\lambda,x)$ will have jumps on $\Gamma$ only.
\paragraph{Discrete spectrum} Choosing $\mathcal H=\ell^2$ and $A = \operatorname{diag}(ik_n^2)$ one construct
in \cite{bib:GenVessel} $S(\lambda,x)$ whose poles are precisely at the points $k_n$.
\paragraph{Fadeyev inverse scattering} Using the ideas of the continuous and the discrete spectrum the Fadeyev inverse
scattering theory \cite{bib:FaddeyevII} was implemented in \cite{bib:GenVessel}.


\bibliographystyle{alpha}
\bibliography{../biblio}

\end{document}